\documentclass
[conference,  10pt,  one column]{IEEEtran}
\usepackage{amsfonts}
\usepackage{amssymb}
\usepackage{graphicx}
\usepackage{amsmath}
\usepackage{epsfig}
\usepackage{epstopdf}
\usepackage{setspace}
\usepackage{color}
\usepackage{cite}
\usepackage{amsmath}
\usepackage{booktabs}
\usepackage[caption=false, font=normalsize, labelfont=sf, textfont=sf]{subfig}

\newcommand{\vv}[1]{\ensuremath{{\vec{#1}}}} 
\newcommand{\sr}[1]{\ensuremath{\vec{x}_{#1}}}
\newcommand{\si}[1]{\ensuremath{\vec{y}_{#1}}}

\newcommand{\Tr}{\text{Tr}}
\renewcommand{\Re}{\text{Re}}
\renewcommand{\Im}{\text{Im}}
\newcommand{\hh}{\boldsymbol{h} }

\newcommand{\A}{\boldsymbol{A}}

\newcommand{\anv} {\ensuremath{\alpha}}
\newcommand{\Tt} {\ensuremath{\tau}}
\newcommand{\av} {\ensuremath{\vec{\alpha}}}
\newcommand{\xv} {\ensuremath{\vec{x}}}
\newcommand{\yv} {\ensuremath{\vec{y}}}
\newcommand{\ev} {\ensuremath{\vec{\epsilon}}}

\newcommand{\aG} {\ensuremath{\alpha}}

\newcommand{\dG} {\ensuremath{\delta}}

\newcommand{\eG} {\ensuremath{\epsilon}}

\newcommand{\gG} {\ensuremath{\gamma}}
\newcommand{\lG} {\ensuremath{\lambda}}

\newcommand{\qG} {\ensuremath{\theta}}

\newcommand{\tG} {\ensuremath{\tau}}

\newcommand{\grad}{\ensuremath{\vv{\nabla}}}

\newcommand{\avp}{\ensuremath{\vec{\alpha}\, '}} 

\begin{document}

\title{Very Low-Complexity Algorithms for Beamforming in Two-Way Relay Systems }

\author{
\authorblockN{Christopher Thron}
\authorblockA{ \
Texas A\&M University, Central Texas\\
Email:thron@ct.tamus.edu
}
\and
\authorblockN{Ahsan Aziz}
\authorblockA{
National Instruments\\
Austin, TX\\
}
}
\maketitle

\begin{abstract}
In this paper,  we present a novel solution for optimal beamforming in a two-way relay (TWR) systems with perfect channel state information. The solution makes use of properties of quadratic surfaces to simplify the solution space of the problem to $\mathbb{R}^4$, and enables the formulation of a differential equation that can be solved numerically to obtain the optimal beamforming matrix.
\end{abstract}
\section{Introduction}
\label{sec:intro}
 Two-way relay (TWR) systems that employ beamforming techniques enable information exchange  with greatly reduced spectral resource requirements  compared to one-way relaying \cite{ANC}. In this paper,  we consider two-way relays with multiple antennas that communicate with two source nodes,  each with one antenna. We also assume that the  channel vectors that determine signal transfer between the relay and the two source nodes is known to the system.  
Existing optimal beamforming algorithms for this system (such as that in \cite{Rui09}) have high computational complexity.
 In this paper,  we present a numerical solution to the optimal beamforming problem which has greatly reduced complexity over previous known solutions. The solution makes use of properties of quadratic surfaces to transform the problem into a differential equation, which can then be expeditiously solved using numerical methods.

The rest of the paper is organized as follows. In Section \ref{sec:system model},  we present the system model and formulate the mathematical optimization problem which specifies the beamforming matrix. In Section \ref{sec: ProbDef}, we show  that the problem can be transformed to an optimization problem with real coefficients, whose solution is a $2 \times 2$ matrix. In Section \ref{sec:RealSolution}, we show that this simplified optimization problem has a solution which is a real matrix.  In Section \ref{sec:Numerical}, we show how the optimization problem can again be transformed into a vector differential equation in $\mathbb{R}^6$, which may be solved numerically using standard methods.

\section{System Model and Formulation of the Optimization Problem}\label{sec:system model}
\begin{figure}
  \centering
  \includegraphics[width=0.35\textwidth]{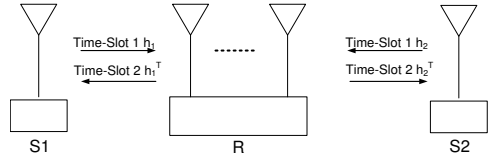}
  \caption{System model}\label{fig: SystemModel}
\end{figure}

We consider a two-way relay system similar to the one introduced in \cite{Rui09},  which consists of the relay node $R$ and two terminal nodes $S1$ and $S2$,  as shown on
Fig. \ref{fig: SystemModel}.

The relay is equipped with $M$ antennas and the terminal nodes are each equipped with a single antenna. Based on the principle
of analog network coding \cite{ANC},  the two terminal nodes exchange information in two consecutive time slots via the help of $R$. In the
first time slot,  terminal nodes $S1$ and $S2$ send messages $s_1$ and $s_2$ with power levels $p_1$ and $p_2$ respectively to $R$,  and the received signal at $R$ is given as
\begin{equation}
\boldsymbol {y}_R=\hh _1\sqrt{p_1}s_1+\hh _2\sqrt{p_2}s_2+\boldsymbol{z}_R, 
\end{equation}
where $\hh _1,  \hh _2 \in \mathcal{C}^{M\times1 }$ are complex channel gains from the terminal nodes $S1$ and $S2$ to the relay
respectively,  $\boldsymbol{z}_R$ is the circularly symmetric complex Gaussian (CSCG) noise with covariance $\sigma_{R}^{2}\boldsymbol{I}$,  and
$E[s_i]=1$,  $i=1, 2$. In the second time slot,  the relay $R$ multiplies a beamforming matrix $\A$ with the received signal $\boldsymbol {y}_R$ and transmits
the resulting vector signal $\A\boldsymbol{y}_R$ to the two terminal nodes. Based on the assumption of channel reciprocity \cite{Zeng11},  the received
signals at $S1$ and $S2$ are given as
\begin{eqnarray}
y_1=\hh _{1}^T\A\hh _1 \sqrt{p_1}s_1 + \hh _{1}^T\A\hh _2 \sqrt{p_2}s_2
+\hh _{1}^T\A\boldsymbol{z}_R+z_1, \label{eq: Source_y1}\\
y_2=\hh _{2}^T\A\hh _2 \sqrt{p_2}s_2 + \hh _{2}^T\A\hh _1 \sqrt{p_1}s_1
+\hh _{2}^T\A\boldsymbol{z}_R+z_2, \label{eq: Source_y2}
\end{eqnarray}
where $z_{1}$ and $z_{2}$ are the CSCG noises at $S1$ and $S2$ with variances $\sigma_{1}^{2}$ and $\sigma_{2}^{2}$,  respectively. In the ideal CSI case as in
\cite{Rui09},  $S1$ and $S2$ can cancel out the self-interference terms $\hh _{1}^T\A\hh _1 \sqrt{p_1}s_1$
and $\hh _{2}^T\A\hh _2 \sqrt{p_2}s_2$ from $y_{1}$ and $y_{2}$,  respectively. 
The corresponding transmit power at the relay $R$ is given by 
\begin{align}
\label{eq:txpwr}
G(\A) \equiv \|\A\hh _1\|^2 p_1+\|\A\hh _2\|^2 p_2+
\Tr[\A^H\A]\sigma_{R}^{2}, 
\end{align}
and the SINRs at node $S_i$ are given by ($i = 1, 2;~k=3-i$)
\begin{align}
\label{eq:snr1}
\text{SINR}_i(\A)= \frac{|{\hh }_{i}^T\A{\hh }_{k}
|^2 p_{k}}
{|\|{\hh }_{i}^T \A\|^2\sigma^2_R+\sigma_i^2}.
\end{align}
Based on these definitions,  the nonrobust optimization problem to minimize the relay power under SINR constraints can be formulated as follows: find (i=1, 2)
\begin{equation}
\label{modOptProblem0}
\A_* = \textrm{arg} \min_{\A}\left[ G(\A)\right] \qquad
 {\rm s.t.} \quad{f_{i}(\A)}\geq \gamma_i \sigma_i^2, 
\end{equation}
where  $\gamma_i$ is the SINR target at $S_i$  and
\begin{equation}
\label{eq:f_i_def}
f_{i}(\A) \equiv  |{\hh }_{i}^T\A{\hh }_{k}|^2p_k  -  |\|{\hh }_{i}^T \A\|^2\sigma^2_R \gamma_i, ,~(k \equiv 3-i).
\end{equation} 

We note that the problem in (\ref{modOptProblem0}) is not convex in general,  because the constraints are not convex functions.

\section{Reduction to  rank 2 problem with real coefficients}\label{sec: ProbDef}
In this section we show how (\ref{modOptProblem}) can be transformed into a much simpler problem with real coefficients. 

It has been shown previously in \cite{Rui09} that $\A_*$ is of complex rank 2.  Specifically,   $\A_*$ can be expressed as
\begin{equation}
\A_* = \sum_{i, j=1}^2 (a_{*})_{ij} \bar{\hh}_i \hh_j^{H}  =  [ \bar{\hh}_1,~\bar{\hh}_2] a_* [\hh_1^{H}~;~ \hh_2^H], 
\end{equation}
where $a_*$ is a complex $2 \times 2$ matrix.  
The objective function condition and constraints in (\ref{modOptProblem0}) can be rewritten in terms of the matrix $a_*$. Note that the coefficients which appear in this simplified version of (\ref{modOptProblem0}) will be complex in general; but it is possible to further simplify the expressions so that all coefficients are real.
After simplification,
the optimization problem becomes:
\medskip
\begin{equation}
\label{modOptProblem}
\anv_* = \textrm{arg} \min_{\anv}\left[ G(\anv)\right] \qquad
 {\rm s.t.} \quad{f_{i}(\anv)}\geq 1, ,~(i=1, 2), 
\end{equation}
where
\begin{equation}\label{eq:power}
G(\anv) \equiv q_1 \| \anv \Tt_1\| ^2 + q_2\| \anv \Tt_2\|^2 + \Tr[\anv^H \anv];
\qquad  f_i(\anv) \equiv c_i|\Tt_i^T\anv\Tt_k|^2 - d_i\|\Tt_i^T\anv\|^2,  \quad (i = 1, 2; k = 3-i),
\end{equation}
where $q_i, c_i$, and $d_i$ are constants ($i=1,2$); and $\tG_i = [1~~\pm r]^T$ where $r$ is a positive real number.
 First we define  $(i = 1, 2; k \equiv 3-i)$
\begin{equation}\label{eq:phidef}
 \Tt_{ii} \equiv \Tt_i \Tt_i^T; \qquad
m  \equiv q_1 \Tt_{11} + q_2 \Tt_{22}  +  I. 
\end{equation}
Next,  for any $2 \times 2$ matrix $A$ we define the operations:
\begin{equation}
\vec{A} \equiv [A_{11}~A_{12}~A_{21}~A_{22}]^T;\qquad \underline{A} \equiv 
\left[ \begin{array}{cc}
A & 0 \\
0 & A  \end{array} \right]; \qquad
\widetilde{A} \equiv 
\left[ \begin{array}{cc}
A_{11}I & A_{21}I \\
A_{12}I & A_{22}I  \end{array} \right].
\end{equation}
Finally we define
\begin{equation}
M \equiv \underline{m}; \qquad T_{ki} \equiv \underline{\tau_{kk}} \widetilde{\tau_{ii}};\qquad  Q_i \equiv c_i T_{ki} -  d_i \widetilde{\tG_{ii}}, 
\end {equation}
where $M,  T_{ki}$,  and $Q_i$ are all real symmetric matrices.
Using this notation,  we have
\begin{equation}\label{eq:power2a}
G(\anv) \equiv  \av^H M \av; \qquad
 f_i(\anv)  \equiv  \av^H Q_i\av,
\end{equation}
where $M$ and $Q$  are real $4\times4$ matrices. With the additional notation
\begin{equation}
\sr{A} \equiv \Re[\vec{A}]; \qquad
\si{A} \equiv \Im[\vec{A}],
\end{equation}
we may rewrite as
\begin{equation}\label{eq:power2}
G(\anv) \equiv  \sr{\anv}^T M\sr{\anv} + \si{\anv}^T M\si{\anv}; \qquad
 f_i(\anv)  \equiv  \sr{\anv}^T Q_i\sr{\anv} + \si{\anv}^T Q_i\si{\anv}.
\end{equation}
 In the following section,  we will show that there always exists an optimal solution $\anv_*$ for (\ref{modOptProblem}) that is also real (so that $(\si{\anv})_* = 0$).

\section{Existence of real optimal solutions }\label{sec:RealSolution}
In this section we show that given a locally-optimal complex feasible solution to (\ref{modOptProblem},   \ref{eq:power2}),   there also exists a  \emph{real} feasible solution  that achieves the \emph{same} power. Since any global optimum is also a local optimum,  it follows that there always exists a globally optimal \emph{real} feasible solution. 

Let us write $\xv \equiv \sr{\anv}$ and $\yv \equiv \si{\anv}$. Then we have
\begin{align}
\label{eq:RealMin2}
G(\xv, \yv) = \xv^T M \xv + \yv^T M \yv; \qquad
f_i(\xv, \yv) =   \xv^T Q_i \xv + \yv^T Q_i \yv, \quad (i=1, 2).
\end{align}
 The following Karush-Kuhn-Tucker (KKT) conditions are satisfied:
\begin{equation}
\label{eq:KKT_simple}
\begin{aligned}
\grad_{\xv}(G(\xv, \yv)) &= \lG_1 \cdot \grad_{\xv}(f_1(\xv, \yv )) + \lG_2 \cdot \grad_{\xv} (f_2(\xv, \yv )) \qquad (\lG_1,  \lG_2 \ge 0);\\
\grad_{\yv} (G(\xv, \yv )) &= \lG_1 \cdot \grad_{\yv} (f_1(\xv, \yv )) + \lG_2 \cdot \grad_{\yv} (f_2(\xv, \yv ))
\end{aligned}
\end{equation}
where
\begin{equation}
\label{eq:KKT_simple_a}
\lG_i(1 - f_i(\xv, \yv ))  = 0,~(i=1, 2).
\end{equation} 
Using (\ref{eq:RealMin2}) we find the explicit KKT conditions are:
\begin{align}
\label{eq:Mxy1}
 M\xv &= \lG_1 Q_1 \xv  + \lG_2 Q_2\xv ; \qquad
M \yv = \lG_1 Q_1 \yv  +   \lG_2 Q_2 \yv,  \quad (\lG_1,  \lG_2 \ge 0).
\end{align}
Substituting (\ref{eq:Mxy1}) into  (\ref{eq:RealMin2}),  we find that the  power achieved at the locally-optimal complex feasible solution is
\begin{equation}
\label{eq:RealMin3}
 \text{Power} = G(x,y) = \lG_1 \left(\xv^T Q_1 \xv   + \yv^T Q_1 \yv  \right) 
+   \lG_2 \left( \xv^T Q_2 \xv   + \yv^T Q_2 \yv   \right).
\end{equation}

Consider first the case where  the constraints ${f_{i}(\anv)}\geq 1$ are both satisfied with equality:
\begin{equation}\label{eq:RealMinConstraint2a}
1 = {f_{i}(\anv)} =    \xv^T Q_i \xv + \yv^T Q_i \yv, \qquad i=1, 2.
\end{equation}
It follows from (\ref{eq:RealMin3}) that
\begin{equation}
\text{Power} =  \lG_1 + \lG_2.
\end{equation}
From (\ref{eq:Mxy1})  we find that the \emph{real} beamforming matrix $\gG_x \xv  + \gG_y \yv $ also satisfies the KKT conditions,  for \emph{any} choice of $\gG_x$ and $\gG_y$:
\begin{align}
  M(\gG_x \xv  + \gG_y \yv ) = &\lG_1 Q_1(\gG_x \xv  + \gG_y \yv ) +   \lG_2  Q_2( \gG_x \xv  + \gG_y \yv ).\label{eq:Mxy3}
\end{align}
Suppose we can find $\gG_x$ and $\gG_y$ so that the SNR constraints are satisfied with equality:
\begin{align}\label{eq:RealMinConstraint3a}
1 &=  (\gG_x \xv  + \gG_y \yv )^T Q_i (\gG_x \xv  + \gG_y \yv )   = \gG_x^2 \cdot  \xv^T Q_i \xv  + \gG_y^2  \cdot \yv^T Q_i \yv  + 2\gG_y \gG_x  \cdot \yv Q_i \xv,  \quad (i=1, 2).
\end{align}
Then the resulting power for this beamforming matrix is also $\lG_1 + \lG_2$,  as above. Thus this real beamforming matrix is feasible,  and achieves the same power as the complex solution. Hence if the complex solution is a global optimum,  it follows that the real solution is a global optimum as well.

It remains to show that it is indeed possible to find $\gG_x,  \gG_y$ that satisfy both constraints in (\ref{eq:RealMinConstraint3a}).
It follows from the fact $\xv^T Q_i \xv  +  \yv^T Q_i \yv = 1$ for $i=1, 2$,  that there are essentially three cases to consider:
\begin{align*}
&\text{(A)}  ,~0 < \xv^T Q_1 \xv , ~\yv^T Q_1 \yv   < 1,  \quad 0 < \xv^T Q_2 \xv ,,~  \yv^T Q_2 \yv   < 1; \\
&\text{(B)}  ,~0 < \xv^T Q_1 \xv ,~ \yv^T Q_1 \yv   < 1,  \quad  \xv^T Q_2 \xv   \ge 1, ,~ \yv^T Q_2 \yv   \le 0; \\
&\text{(C)}  ,~\xv^T Q_1 \xv   \ge 1,~ \yv^T Q_1 \yv   \le 0,  \quad   \xv^T Q_2 \xv   \le 0,,~  \yv^T Q_2 \yv   \ge 1.
\end{align*}
All other cases can be reduced to one of these cases by exploiting the symmetry between $\xv $ and $\yv $,  and between $\phi_1$ and $\phi_2$.  Notice that the case $\xv^T Q_1 \xv,~ \xv^T Q_2 \xv  > 1$ is impossible, since then  we would have  $\yv^T Q_1 \yv,~\yv^T Q_2 \yv < 0$ so that $\av = \xv $ satisfies both constraints and is a feasible solution with lower power than $\av = \xv  + j\yv $. 

In case (A),   (\ref{eq:Mxy3}) yields ellipses in the $(\gG_x, \gG_y)$ plane for $i=1,2$.  The positive $\gG_x$-intercepts for the two constraints are $(\xv^T Q_1 \xv  )^{-1/2}$ and $(\xv^T Q_2 \xv  )^{-1/2}$ respectively; while the positive $\gG_y$-intercepts for the two constraints are $(\yv^T Q_1 \yv  )^{-1/2}$ and $(\yv^T Q_2 \yv  )^{-1/2}$ respectively. However,  from
(\ref{eq:RealMinConstraint2a}) we have $\yv^T Q_1 \yv   = 1 - \xv^T Q_1 \xv  $. Hence the order of positive $\gG_y$-intercepts for the two constraints is the reverse of the order for positive $\gG_x$-intercepts. It follows that the two constraint ellipses must cross somewhere in the first quadrant. At the crossing point,  both constraints are satisfied with equality. 

In case (B),  the first constraint corresponds to an ellipse and the second to a hyperbola  in the $(\gG_x, \gG_y)$ plane.. The  positive $\gG_x$ intercept for the elliptical constraint is $(\xv^T Q_1 \xv  )^{-1/2} > 1$,  while the  positive $\gG_x$ intercept for the hyperbolic constraint is $(\xv^T Q_2 \xv  )^{-1/2} \le 1$. Since the ellipse encloses at least one point on the hyperbolic constraint and the hyperbolic constraint is unbounded,  it follows that the elliptical and hyperbolic constraints must intersect,  so there must be at least one point where both constraints are satisfied with equality.

Case (C) can actually be reduced to Case (A) or Case (B).  Note that if $\av$ is a solution, then $\av_{(\qG)} \equiv e^{i\qG}\av$ is also a solution. We have:
\begin{equation}
\xv_{(\qG)} = \cos \qG \cdot \xv - \sin \qG \cdot \yv; \qquad  \yv_{(\qG)} = \sin \qG \cdot \xv + \cos \qG \cdot \yv.
\end{equation}
 It follows that 
\begin{align}
\xv_{\anv_{(\qG)}}^T Q_1 \xv_{\anv_{(\qG)}} &= (\cos \qG \cdot \xv - \sin \qG \cdot \yv)^T Q_1 (\cos \qG \cdot \xv - \sin \qG \cdot \yv) \notag  \\
&= \cos^2  \qG \cdot \xv^T Q_1  \xv + (1 - \cos^2 \qG) \cdot \yv^T Q_1 \yv + \sin 2\qG \cdot \xv Q_1 \yv.
\end{align}
Clearly $\xv_{\anv_{(\qG)}}^T Q_1 \xv_{\anv_{(\qG)}}$ is a continuous function of $\qG$.  When $\qG = 0$, in case (C) we have   $\xv_{\anv^{(0)}}^T Q_1 \xv_{\anv^{(0)}} 
= \xv^T Q_1 \xv \ge 1$  However, when  $\qG = \pi/2$ we have $\xv_{\anv^{(\pi/2)}}^T Q_1 \xv_{\anv^{(\pi/2)}} = (1 -  \xv^T Q_1 \xv)  \le 0$. By continuity, there must be a value of $\qG$ such that $0 < \xv_{\anv_{(\qG)}}^T Q_1 \xv_{\anv_{(\qG)}} < 1$: and case (A) or (B) applies in this situation. This completes the argument in the case where both constraints in (\ref{eq:RealMinConstraint3a}) are satisfied with equality.

It is also possible that  only one of the constraint conditions in (\ref{modOptProblem}) is satisfied with equality. In this case,  then similar  arguments can be used in cases (A) and (B). 
(The above argument for case (C) also holds if only one of the constraints holds with equality.)  In case (A)
we may suppose that the $i=1$  constraint holds with equality,  while the  $i=2$ constraint holds with strict  inequality.  It follows that either  $\xv^T Q_1 \xv  <  \xv^T Q_2 \xv  $ or  $\yv^T Q_1 \yv   <  \yv^T Q_2 \yv  $. Without loss of generality,  we may suppose that 
$\xv^T Q_1 \xv  <  \xv^T Q_2 \xv  $. In this case,  then
using $ \gG_x  = (\xv^T Q_1 \xv  )^{-1/2}$ and $\gG_y = 0$ yields a real solution that also satisfies the $i=1$ constraint with equality and the $i=2$ constraint with strict inequality,  and has the same power. Case (B) must be divided into two cases. In the case where the elliptic constraint (which we may assume corresponds to $i=1$) holds with equality and the hyperbolic constraint (corresponding to $i=2$) holds with strict inequality,  then similar arguments show that   $ \gG_x = (\xv^T Q_1 \xv  )^{-1/2}$ and $\gG_y  = 0$ yields a real solution with the same power that satisfies both constraints. If the hyperbolic constraint holds with equality and the elliptical constraint with strict inequality,   then since the hyperbolic constraint is unbounded it is always possible to find  $\gG_x$ and $\gG_y$ such that the elliptic constraint is satisfied. 
 
It is not possible  for both constraint conditions to hold with strict inequality,  since then
(\ref{eq:KKT_simple_a}) gives $\lG_1 = \lG_2 = 0$,  so (\ref{eq:KKT_simple}) implies that $\grad_{\xv} G(\xv, \yv ) = \grad_{\yv} G(\xv, \yv ) = 0$,  which in turn implies that $\xv =\yv =0$ since $G(\xv, \yv )$ is positive definite. 

In summary,  we have shown that there always exists a  real optimum solution. This reduces the complexity of the problem by a factor of more than 2, since a complex addition requires 2 real additions, while a complex multiplication requires 4 real multiplications.

\section{Numerical solution to the reduced problem}\label{sec:Numerical}

\subsection{Exact solution for case $d_i=0$}
An exact solution to (\ref{modOptProblem}) is possible in the case where $d_i=0,~(i=1,2)$. According to the results of the previous section,  we may assume that the solution $\anv$ is real. The constraint inequalities $f_i(\anv) \ge 1$ become (from (\ref{eq:power})):
\[
|\Tt_i^T\anv\Tt_k|^2  \ge 1/c_i, \qquad (i=1, 2;~k=3-i).
\]
Writing out these constraints in terms of matrix components (and replacing inequality with equality) gives:
\begin{align*}
\anv_{11} + r(\anv_{21} - a_{12}) - r^2 \anv_{22} =   c_1^{-1/2}; \qquad
\anv_{11} + r(\anv_{12} - \anv_{21}) -  r^2 \anv_{22} = \pm c_2^{-1/2}, 
\end{align*}
where without loss of generality we have chosen the positive sign in the first equation  since the optimal beamforming matrix is arbitrary up to an overall minus sign. In vector notation, this becomes
\begin{equation}
[1,~-r ,~ r ,~- r^2] \av = c_1^{-1/2};\qquad  [1,~r ,,~-r ,~- r^2] \av = c_2^{-1/2}.
\end{equation}
These constraints correspond to a pair of parallel hyperplanes in $\mathbb{R}^4$, which intersect in four two-dimensional planes  as long as all hyperplanes are not parallel (which can only occur if $\hh_1 = \hh_2$). 
The sets in $\mathbb{R}^4$ that correspond to constant power are a concentric family of ellipsoids centered at the origin. The minimum-power solution corresponds to the smallest ellipsoid that touches at least one of these planes. This geometrical argument shows that the optimal solution will satisfy both constraints with equality. 

The set in $\mathbb{R}^4$ which satisfies these constraints with equality is a plane which is given parametrically by ($z_1,  z_2$ are arbitrary real parameters)
\begin{equation}
\av^T = 
\left[ 
\frac{c_1^{-1/2} \pm c_2^{-1/2}}{2} ,~\frac{c_1^{-1/2} \mp c_2^{-1/2}}{2r},~0,~0 \right] 
+ z_1 [ r^2 ,~ 0 ,~0 ,~1 ]
+ z_2 [ 0 ,~ 1 ,~1 ,~ 0 ].  
\end{equation}
Power is minimized when the derivatives of $G(\av) \equiv\av^T  M \av$ with respect to $z_1$ and $z_2$ are equal to zero.  This gives the two conditions
\begin{equation}
[ r^2 ,~ 0 ,~0 ,~1 ] M \av = [ 0 ,~ 1 ,~1 ,~ 0 ] M\av = 0. 
\end{equation}
These equations gives two different solutions for $z_1,  z_2$ corresponding to the two different sign choices for $\pm c_2^{-1/2}$.

In summary,  we have computed two real candidate optimal solutions  for the case $d_i=0$. The overall optimal solution will be the candidate which has the lowest power.

\subsection{Numerical solution for $d_i \neq 0$}
In order to obtain solutions for $d_i \neq 0$,  we assume that the equations (\ref{modOptProblem}) have been solved for the case where $d_i$ is replaced by $w d_i$,  where $w$ is a parameter between 0 and 1. We may then  use the solution for $wd_i$ to find the solution for $(w+\dG) d_i$,  where $\dG$ is an incremental change in the value of $w$.  In this way,  we may obtain a differential equation for the solution with arbitrary $wd_i$,  using one of the solutions for $w = 0$  as initial conditions. Plugging in $w=1$,  we obtain a solution corresponding to the given $d_i$.  We may use as initial conditions either of the two solutions corresponding to $d_i=0$ described in the preceding section.

This process of increasing $w$ from 0 to 1 has a geometrical interpretation. Consider $\av$ as an element of $\mathbb{R}^4$.  Then when $w=0$,  each constraint  $f_i(\av) = 1~(i=1,2)$ corresponds to a pair of parallel 3-dimensional hyperplanes. Each hyperplane for $f_1(\av) = 1$ intersects each hyperplane for $f_2(\av) = 1$ in a 2-dimensional plane lying in $\mathbb{R}^4$. There are thus four 2-dimensional planes in $\mathbb{R}^4$ where both constraints are satisfied with equality. Two of these planes are the negatives of the other two,  so we need only consider two of these planes. At the same time,  the constant-power surfaces correspond to a family of concentric,  disjoint ellipsoids in $\mathbb{R}^4$ centered at the origin. The smallest of these ellipsoids that intersects at least one of the two 2-dimensional planes will be tangent at a single point because of the strict convexity of the ellipsoids.This single point is the optimal beamforming matrix in the case where $w = 0$.  If we now  let $w$ increase,  each constraint hyperplanes ``bends'' and  become one sheet of a  hyperboloid. the intersection of each sheet for $i=1$  with each $i=2$ sheet is either empty or  a  2-dimensional hyperboloidal surface. Thus the set of intersections consists of at most 4 2-dimensional hyperboloidal surfaces. Since the negative of each surface of intersection is also a surface of intersection, there are at most two surfaces that need to be considered. Because of convexity properties of these 2-dimensional surfaces,  the smallest ellipse that intersects at least one of these surfaces will intersect at a single point,  which is the optimal beamforming solution. 

Let $\av$ be the optimal beamforming matrix for $wd_i$,  and let $\av + \ev$ be the perturbed solution corresponding to $(w+\dG)d_i$. We now derive a first-order expression for $\ev$,   which will lead to a differential equation for $\av$ as a function of $w$ on the interval $0 \le w \le 1$.  
For ease of notation, we define
\begin{equation}
Q_i^{(w)} \equiv  c_i T_{ki} -  w d_i \widetilde{\tG_{ii}}, \quad (i=1,2).
\end{equation}
Thus $Q_i$ in the preceding discussion corresponds to $Q_i^{(1)}$, while the $d_i=0,~(i=1,2)$ case corresponds to $Q_i^{(0)}$. 

The constraint equations corresponding to $(w+\dG)d_i$  become:
\begin{align*}
 (\av + \ev)^T  Q_i^{(w+\dG)}  (\av + \ev)  = 1 ,~~ (i=1, 2), 
\end{align*}
where we have presented the constraints with equality because our above argument establishes that both constraints will be met with equality. The equation is satisfied to zeroth order by assumption,  and to first order we have
\begin{align}\label{eq:fo_constraint}
2 \av^T Q_i^{(w)} \ev  =    d_i \av^T  \av \dG   , \qquad (i=1, 2).
\end{align}
The maximization (KKT) condition is
\begin{align*}
\grad_{\eG} G(\anv + \eG) = \lG_1'  \grad_{\eG} f_1(\anv + \eG ) + \lG_2'  \grad_{\eG} f_2(\anv + \eG), 
\end{align*}
which written out more explicitly is

\begin{align}\label{eq:power3}
&\grad_{\eG}  (\av + \ev)^T M (\av + \ev)  =
\sum_{i=1}^2 \lG_i'  \grad_{\ev}  (\av + \ev)^T Q_i^{(w+\dG)} (\av + \ev)  ] \quad ( k = 3-i).
\end{align}
In order to write this perturbatively,  we define
\[ \eta_i = \lG_i' - \lG_i \qquad (i=1, 2). \]
Then to first order,  (\ref{eq:power3}) becomes (note all zeroeth-order terms cancel)
\begin{equation}\label{eq:power4}
 M \ev   = \sum_{i=1}^{2}   \lG_i  Q_i^{(w)}\ev -   d_i \widetilde{\tG_{ii}} \av \dG +  Q_i^{(w)}\av \eta_i  , 
\end{equation}
which can be rearranged to give
\begin{equation}\label{eq:power4}
  \left(\sum_{i=1}^{2}   \lG_i  Q_i^{(w)} - M\right)\ev   +  Q_i^{(w)}\av \eta_i  = d_i \widetilde{\tG_{ii}} \av \dG  
\end{equation}
We may now replace $\ev / \dG$ with $\frac{d \av}{d w}$ and $\eta_i/ \dG$ with $\frac{d  \lG_i}{d w}$
in equations (\ref{eq:fo_constraint}) and (\ref{eq:power4})
 to obtain a system of six ordinary differential equations for the four entries of $\av$ plus the two Lagrange multipliers $\lG_1, \lG_2$.

\begin{align}\label{eq:sysODE}
  \left(\sum_{i=1}^{2}   \lG_i  Q_i^{(w)} - M\right)\frac{d \av}{dw}   +  Q_i^{(w)}\av \frac{d \lG_i}{dw}  = d_i \widetilde{\tG_{ii}} \av ; \\
\av^T Q_i^{(w)}  \frac{d \av}{dw} = \frac{d_i}{2} \av^T \widetilde{\tG_{ii}} \av.
\end{align}
We may rewrite this in more conventional form as a vector ODE. Letting $\avp$ and $\lG'_i$ denote  $\frac{d \av}{d w}$ and $\frac{d  \lG_i}{d w}$ respectively, 
we may rewrite the system as:
\begin{align}\label{eq:sysMx}
\left[ \begin{array}{ccc}
 \sum_{i=1}^{2}  \lG_i Q_i^{(w)}  -M &Q_1^{(w)} \vec{\aG} & Q_2^{(w)} \vec{\aG}\\
 (Q_1^{(w)} \vec{\aG})^T  & 0 & 0 \\
( Q_2^{(w)} \vec{\aG})^T  & 0 & 0 
 \end{array} \right]
\left[ \begin{array}{c}
  \avp \\ \lG'_1 \\ \lG'_2
 \end{array} \right]
=
\left[ \begin{array}{c}
\sum_{i=1}^{2} \lG_i  d_i   (\widetilde{\tG_{ii}} \vec{\aG}) \\
 \frac{d_1}{2} \vec{\aG}^T(\widetilde{\tG_{11}} \vec{\aG}) \\
 \frac{d_2}{2} \vec{\aG}^T(\widetilde{\tG_{22}} \vec{\aG})
 \end{array} \right]
\end{align}

We may use the Runge-Kutta method  to solve this on the interval $0 \le w \le 1$.  The initial conditions are obtained from the $w=0$ solution obtained above.

\bibliographystyle{IEEEtran}
\bibliography{ref}
\end{document}